\documentclass[oneside,french,english]{amsart}
\usepackage[T1]{fontenc}
\usepackage[utf8]{inputenc}
\usepackage[a4paper]{geometry}
\usepackage{babel}
\makeatletter
\addto\extrasfrench{%
   \providecommand{\fg}{\ifdim\lastskip>\z@\unskip\fi~\frqq}%
}

\makeatother
\usepackage{url}
\usepackage{amsbsy}
\usepackage{amstext}
\usepackage{amsthm}
\usepackage{amssymb}
\usepackage{graphicx}
\usepackage[unicode=true]
 {hyperref}

\makeatletter
\numberwithin{equation}{section}
\numberwithin{figure}{section}

\usepackage{tikz}
\usepackage{tikz-3dplot}
\usetikzlibrary{patterns}
\usetikzlibrary{matrix}
\usetikzlibrary{decorations.pathmorphing}
\usetikzlibrary{decorations.pathreplacing}
\usetikzlibrary{decorations.markings}
\usetikzlibrary{decorations.shapes}
\usetikzlibrary{arrows}
\usepackage{pgfplots}
\pgfplotsset{compat=newest}

\makeatother

\begin{document}

\title{Task-based parallelization of an implicit kinetic scheme}

\author{Jayesh Badwaik, Matthieu Boileau, David Coulette, Emmanuel Franck,\\Philippe
Helluy, Laura Mendoza, Herbert Oberlin}

\address{\href{mailto:helluy@unistra.fr}{helluy@unistra.fr}}
\begin{abstract}
In this paper we present and implement the Palindromic Discontinuous
Galerkin (PDG) method in dimensions higher than one. The method has
already been exposed and tested in \cite{coulette2016palindromic}
in the one-dimensional context. The PDG method is a general implicit
high order method for approximating systems of conservation laws.
It relies on a kinetic interpretation of the conservation laws containing
stiff relaxation terms. The kinetic system is approximated with an
asymptotic-preserving high order DG method. We describe the parallel
implementation of the method, based on the StarPU runtime library.
Then we apply it on preliminary test cases.
\end{abstract}

\maketitle

\section{Introduction}

\global\long\def\vf{\mathbf{{f}}}
\global\long\def\vg{\mathbf{{g}}}
\global\long\def\vh{\mathbf{{h}}}
\global\long\def\vs{\mathbf{{s}}}
\global\long\def\vV{\mathbf{V}}
\global\long\def\vx{\mathbf{x}}
\global\long\def\v#1{\mathbf{#1}}
\global\long\def\vq{\mathbf{q}}
\global\long\def\vw{\mathbf{w}}
\global\long\def\ddim{D}
\global\long\def\dorder{d}
\global\long\def\normal{\mathbf{n}}
\global\long\def\flux{\mathbf{q}}
\global\long\def\Vu{\mathbf{u}}
\global\long\def\Vvi{\mathbf{v}_{i}}
In this work we consider the time discretization of compressible fluid
models that appear in gas dynamics, biology, astrophysics or plasma
physics for tokamaks. These models can be unified in the following
form 
\begin{equation}
\partial_{t}\vw+\sum_{k=1}^{\ddim}\partial_{k}\flux^{k}(\vw)=\vs,\label{eq:hyper}
\end{equation}
where $\vw:\mathbb{R}^{\ddim}\times[0,T_{\max}]\longrightarrow\mathbb{R}^{m}$
is the vector of conservative variables, $\flux^{k}(\vw):\mathbb{R}^{m}\longrightarrow\mathbb{R}^{m}$
is the flux and $\vs:\mathbb{R}^{\ddim}\times\mathbb{R\times R}^{m}\longrightarrow\mathbb{R}^{m}$
is a source term. $\ddim$ represents the physical space dimension
and $m$ the number of unknowns.

In many physical applications such as MHD flows, low Mach Euler equations,
Shallow-Water with sedimentation, the model presents several time
scales associated to the propagation of different waves. When the
time scale of fast phenomena, which constrains the explicit CFL condition,
is very small compared to the time scale of the most relevant phenomena,
it becomes necessary to switch to implicit schemes. However classical
implicit schemes are very costly in 2D or 3D because they require
the resolution of linear or non-linear systems at each time step.
In addition, the matrices associated with the hyperbolic systems are
generally ill-conditioned.

In this paper, we propose to follow another approach for avoiding
the resolution of complicated linear systems. Instead of solving the
full fluid model (\ref{eq:hyper}) directly, we replace it by a simpler
kinetic interpretation made of a set of transport equations coupled
through a stiff relaxation term \cite{aregba2000discrete,bouchut2000kinetic,graille2014approximation}.
See also \cite{coulette2016palindromic} and included references.
The kinetic system is then solved by a splitting method where the
transport and relaxation stages are treated separately. The method
is then well adapted to parallel optimizations. The method is already
presented in \cite{coulette2016palindromic} in the one-dimensional
case. Here we present its implementation in higher dimensions. We
particularly focus our presentation of the massive parallelization
of the method with the StarPU runtime system \cite{AugAumFurNamThi2012EuroMPI}.

The outlines are as follows.

First we recall that it is possible to provide a general kinetic interpretation
of any system of conservation laws. The interest of this representation
is that the complicated non-linear system is replaced by a (larger)
set of scalar linear transport equations that are much easier to solve.
The transport equations are coupled through a non-linear source term
that is fully local in space.

Then, we detail the approximation which allows to solve the transport
equation in an efficient way. We adopt a Discontinuous Galerkin (DG)
method based on an upwind numerical flux and Gauss-Lobatto quadrature
points. Thanks to the upwind flux, the matrix of the discretized transport
operator has a block-triangular structure.

The main part of this work is devoted to the task parallelization
based on the StarPU library for treating the transport and relaxation
steps efficiently. We use the MPI version of StarPU, which allows
to address clusters of multicore processors. We also describe the
domain decomposition and the macrocell approach that we have used
to achieve better performance. 

Finally, we give some numerical results for measuring the efficiency
of the parallelism and the ability of the method to compute realistic
problems.

\section{Kinetic model}

We consider the following kinetic equation
\begin{equation}
\partial_{t}\vf+\sum_{k=1}^{\ddim}\vV^{k}\partial_{k}\vf=\frac{1}{\tau}(\vf^{eq}(\vf)-\vf)+\vg.\label{eq:kin_bgk}
\end{equation}
The unknown is a vectorial distribution function $\vf(\vx,t)\in\mathbb{R}^{n_{v}}$
depending on the space variable $\vx=(x^{1}\ldots x^{\ddim})\in\mathbb{R}^{\ddim}$
and time $t\in\mathbb{R}$. $\vg(\v x,t,\v f)$ is a vectorial source
term, possibly depending on space, time and $\vf$. The partial derivatives
are noted
\[
\partial_{t}=\frac{\partial}{\partial t},\quad\partial_{k}=\frac{\partial}{\partial x^{k}}.
\]
The relaxation time $\tau$ is a small positive constant. The constant
matrices $\vV^{k}$, $1\leq k\leq\ddim$ are diagonal
\[
\vV^{k}=\left(\begin{array}{cccc}
v_{1}^{k}\\
 & v_{2}^{k}\\
 &  & \ddots\\
 &  &  & v_{n_{v}}^{k}
\end{array}\right)
\]

In other words, (\ref{eq:kin_bgk}) is a set of $n_{v}$ transport
equations at constant velocities $\v v_{i}=(v_{i}^{1},\ldots,v_{i}^{D})$,
coupled through a stiff BGK relaxation, and with an optional additional
source term. We denote by $\v{V\cdot\boldsymbol{\partial}}=\sum_{k=1}^{\ddim}\v V^{k}\partial_{k}$
the transport operator, and by $\v{N\vf}=(\vf^{eq}(\vf)-\vf)/\tau$
the BGK relaxation term (also called the ``collision'' term).

Generally, this kinetic model represents an underlying hyperbolic
system of conservation laws. The macroscopic conservative variables
$\v w(\v x,t)\in\mathbb{R}^{m}$ are obtained through a linear transformation
\begin{equation}
\v w=\v P\vf,\label{eq:micro_to_macro}
\end{equation}
where $\v P$ is a $m\times n_{v}$ matrix. Generally the number of
conservative variables is smaller than the number of kinetic data:
$m<n_{v}$. The equilibrium (or ``Maxwellian'') distribution $\vf^{eq}(\vf)$
is such that
\[
\v{Pf}=\v P\vf^{eq}(\vf),
\]
and 
\begin{equation}
\vw=\v P\vf_{1}=\v P\v{f_{2}}\Rightarrow\vf^{eq}(\vf_{1})=\vf^{eq}(\v{f_{2}}),\label{eq:rank_f}
\end{equation}
which states that the equilibrium actually depends only on the macroscopic
data $\vw$. We could have used the notation $\vf^{eq}=\vf^{eq}(\vw)=\vf^{eq}(\mathbf{P}\vf)$,
but we have decided to respect a well-established tradition.

When $\tau\to0$, the kinetic equations provide an approximation of
the system of conservation laws
\begin{equation}
\partial_{t}\vw+\sum_{k=1}^{\ddim}\partial_{k}\flux^{k}(\vw)=\vs,\label{eq:macro_limit_system}
\end{equation}
where the flux is given by

\[
\flux^{k}(\vw)=\v P\vV^{k}\vf^{eq}(\vf).
\]

The flux is indeed a function of $\vw$ only because of (\ref{eq:rank_f}).

Similarly the source term is given by

\begin{equation}
\vs(\vx,t,\vw)=\mathbf{P}\vg(\vx,t,\vf^{eq})
\end{equation}

System (\ref{eq:kin_bgk}) has to be supplemented with conditions
at the boundary $\partial\Omega$ of the computational domain $\Omega$.
We denote by $\normal=(n_{1}\ldots n_{\ddim})$ the outward normal
vector on $\partial\Omega.$ For simplicity, we shall only consider
very simple imposed and time-independent boundary conditions $\vf^{b}$.
We note
\[
\vV\cdot\normal=\sum_{k=1}^{\ddim}\vV^{k}n_{k},\quad\vV\cdot\normal^{+}=\max(\vV\cdot\normal,0),\quad\vV\cdot\normal^{-}=\min(\vV\cdot\normal,0).
\]
A natural boundary condition, which is compatible with the transport
operator $\v{V\cdot\boldsymbol{\partial}}$, is 
\begin{equation}
\vV\cdot\normal^{-}\vf(\vx,t)=\vV\cdot\normal^{-}\vf^{b}(\vx),\quad\vx\in\partial\Omega.\label{eq:boundary_conditions}
\end{equation}
It states that for a given velocity $\v v_{i}$, the corresponding
boundary data $f_{i}^{b}$ is used only at the inflow part of the
boundary.

Let us point out that the programming optimization that we propose
in this paper rely in an essential way on the nature of the boundary
condition (\ref{eq:boundary_conditions}). For other boundary conditions,
such as periodic or wall conditions, additional investigations are
still needed.

\section{Numerical method}

\subsection{Discontinuous Galerkin approximation}

For solving (\ref{eq:kin_bgk}) we shall treat the transport operator
$\vV\cdot\boldsymbol{\partial}$ and the collision operator $\v N$
efficiently, thanks to a splitting approach. This allows to achieve
a better parallelism. Let us start with the description of the transport
solver.

For a simple exposition, we only consider one single scalar transport
equation for $f(\vx,t)\in\mathbb{R}$ at constant velocity $\mathbf{v}$
\begin{equation}
\partial_{t}f+\mathbf{v}\cdot\nabla f=0.\label{eq:single_transport}
\end{equation}
The general vectorial case is easily deduced.

We consider a mesh $\mathcal{M}$ of $\Omega$ made of open sets,
called ``cells'', $\mathcal{M}=\left\{ L_{i},\,i=1\ldots N_{c}\right\} $.
In the most general setting, the cells satisfy
\begin{enumerate}
\item $L_{i}\cap L_{j}=\emptyset$, if $i\neq j$;
\item $\overline{\cup_{i}L_{i}}=\overline{\Omega}.$
\end{enumerate}
In each cell $L\in\mathcal{M}$ we consider a basis of functions $(\varphi_{L,i}(\vx))_{i=0\ldots N_{\dorder}-1}$
constructed from polynomials of order $\dorder$. We denote by $h$
the maximal diameter of the cells. With an abuse of notation we still
denote by $f$ the approximation of $f$, defined by
\[
f(\vx,t)=\sum_{j=0}^{N_{\dorder}-1}f_{L,j}(t)\varphi_{L,j}(\vx),\quad\vx\in L.
\]
The DG formulation then reads: find the $f_{L,j}$'s such that for
all cell $L$ and all test function $\varphi_{L,i}$
\begin{equation}
\int_{L}\partial_{t}f\varphi_{L,i}-\int_{L}f\mathbf{v}\cdot\nabla\varphi_{L,i}+\int_{\partial L}(\mathbf{v}\cdot\normal^{+}f_{L}+\mathbf{v}\cdot\normal^{-}f_{R})\varphi_{L,i}=0.\label{eq:dg_var}
\end{equation}
In this formula (see Figure \ref{eq:convention_downwind}):
\begin{itemize}
\item $R$ denotes the neighboring cell to $L$ along its boundary $\partial L\cap\partial R$,
or the exterior of $\Omega$ on $\partial L\cap\partial\Omega$.
\item $\normal=\normal_{LR}$ is the unit normal vector on $\partial L$
oriented from $L$ to $R$.
\item $f_{R}$ denotes the value of $f$ in the neighboring cell $R$ on
$\partial L\cap\partial R$.
\item If $L$ is a boundary cell, one may have to use the boundary values
instead: $f_{R}=f^{b}$ on $\partial L\cap\partial\Omega$.
\item $\mathbf{v}\cdot\normal^{+}f_{L}+\mathbf{v}\cdot\normal^{-}f_{R}$
is the standard upwind numerical flux encountered in many finite volume
or DG methods.
\end{itemize}
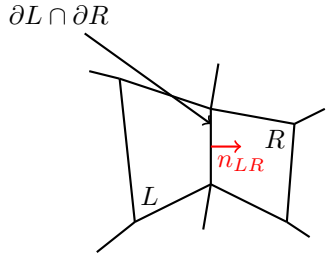
\begin{figure}
\begin{centering}
\begin{center}
\begin{tikzpicture}[scale=1]
\draw[thick] (0,1) -- (0,0); \draw[->, thick, color=red] (0,0.5) -- (0.4,0.5); \node[below, color=red] at (0.4,0.5) {$n_{LR}$}; \node[above] (n1) at (-2,2) {$\partial L\cap\partial R$}; \draw[->, thick] (n1) -- (0,0.8);
\draw[thick] (-1,-0.5) -- (0,0); \draw[thick] (-1.2,1.4) -- (0,1); \draw[thick] (-1.2,1.4) -- (-1,-0.5); \node[above right] at (-1.05,-0.4) {$L$}; 
\draw[thick] (1,-0.5) -- (0,0); \draw[thick] (1.1,0.8) -- (0,1); \draw[thick] (1.1,0.8) -- (1,-0.5); \node[below left] at (1.1,0.85) {$R$};
\draw[thick] (0,0) -- (-0.1,-0.6); \draw[thick] (0,1) -- (0.1,1.6); \draw[thick] (-1,-0.5) -- (-1.5,-0.9); \draw[thick] (-1.2,1.4) -- (-1.6,1.5); \draw[thick] (1,-0.5) -- (1.3,-0.7); \draw[thick] (1.1,0.8) -- (1.5,1);  \end{tikzpicture} 
\par\end{center}
\par\end{centering}
\caption{\label{fig:cell_convention}Convention for the $L$ and $R$ cells
orientation.}
\end{figure}

In our application, we consider hexahedral cells. We have a reference
cell 
\[
\hat{L}=]-1,1[^{\ddim}
\]
\global\long\def\jacob{\boldsymbol{\tau}}
and a smooth transformation $\vx=\jacob_{L}(\hat{\vx})$, $\hat{\vx}\in\hat{L}$,
that maps $\hat{L}$ on $L$
\[
\jacob_{L}(\hat{L})=L.
\]
We assume that $\jacob_{L}$ is invertible and we denote by $\jacob_{L}'$
its (invertible) Jacobian matrix. We also assume that $\jacob_{L}$
is a direct transformation\textcompwordmark{}
\[
\det\jacob_{L}'>0.
\]

In our implementation $\jacob_{L}$ is a quadratic map based on hexahedral
curved ``H20'' finite elements with 20 nodes. The mesh of H20 finite
elements is generated by \texttt{gmsh} \cite{geuzaine2009gmsh}. 

On the reference cell, we consider the Gauss-Lobatto points $(\hat{\vx}_{i})_{i=0\ldots N_{\dorder}-1}$,
$N_{\dorder}=(\dorder+1)^{\ddim}$ and associated weights $(\omega_{i})_{i=0\ldots N_{\dorder}-1}$.
They are obtained by tensor products of the $(\dorder+1)$ one-dimensional
Gauss-Lobatto (GL) points on $]-1,1[$. The reference GL points and
weights are then mapped to the physical GL points of cell $L$ by
\begin{equation}
\vx_{L,i}=\jacob_{L}(\hat{\vx}_{i}),\quad\omega_{L,i}=\omega_{i}\det\jacob_{L}'(\hat{\vx}_{i})>0.\label{eq:map_GL}
\end{equation}
In addition, the six faces of the reference hexahedral cell are denoted
by $F_{\epsilon}$, $\epsilon=1\ldots6$ and the corresponding outward
normal vectors are denoted by $\hat{\normal}_{\epsilon}$. A big advantage
of choosing the GL points is that the volume and the faces share the
same quadrature points. A special attention is necessary for defining
the face quadrature weights. If a GL point $\hat{\vx}_{i}\in F_{\epsilon}$,
we denote by $\mu_{i}^{\epsilon}$ the corresponding quadrature weight
on face $F_{\epsilon}$. We also use the convention that $\mu_{i}^{\epsilon}=0$
if $\hat{\vx}_{i}$ does not belong to face $F_{\epsilon}$. A given
GL point $\hat{\vx}_{i}$ can belong to several faces when it is on
an edge or in a corner of $\hat{L}$. Because of symmetry, we observe
that if $\mu_{i}^{\epsilon}\neq0$, then the weight $\mu_{i}^{\epsilon}$
does not depend on $\epsilon$.

We then consider basis functions $\hat{\varphi_{i}}$ on the reference
cell: they are the Lagrange polynomials associated to the Gauss-Lobatto
point and thus satisfy the interpolation property
\[
\hat{\varphi}_{i}(\hat{\vx}_{j})=\delta_{ij}.
\]
The basis functions on cell $L$ are then defined according to the
formula
\[
\varphi_{L,i}(\vx)=\hat{\varphi}_{i}(\jacob_{L}^{-1}(\vx)).
\]
In this way, they also satisfy the interpolation property
\begin{equation}
\varphi_{L,i}(\vx_{L,j})=\delta_{ij}.\label{eq:interp_property}
\end{equation}
In this paper, we only consider conformal meshes: the GL points on
cell $L$ are supposed to match the GL points of cell $R$ on their
common face. Dealing with non-matching cells is the object of a forthcoming
work.

Let $L$ and $R$ be two neighboring cells. Let $\vx_{L,j}$ be a
GL point in cell $L$ that is also on the common face between $L$
and $R$. In the case of conformal meshes, it is possible to define
the index $j'$ such that
\[
\vx_{L,j}=\vx_{R,j'}.
\]

Applying a numerical integration to (\ref{eq:dg_var}), using (\ref{eq:map_GL})
and the interpolation property (\ref{eq:interp_property}), we finally
obtain
\begin{multline}
\partial_{t}f_{L,i}\omega_{L,i}-\sum_{j=0}^{N_{\dorder}-1}\mathbf{v}\cdot\nabla\varphi_{L,i}(\vx_{L,j})f_{L,j}\omega_{L,j}+\\
\sum_{\epsilon=1}^{6}\mu_{i}^{\epsilon}\left(\mathbf{v}\cdot\normal_{\epsilon}(\vx_{L,i})^{+}f_{L,i}+\mathbf{v}\cdot\normal_{\epsilon}(\vx_{L,i})^{-}f_{R,i'}\right)=0.\label{eq:dg_lobatto}
\end{multline}
We have to detail how the gradients and normal vectors are computed
in the above formula. Let $\v A$ be a square matrix. We recall that
the cofactor matrix of $\v A$ is defined by
\begin{equation}
\text{co}(\v A)=\det(\v A)\left(\v A^{-1}\right)^{T}.\label{eq:co_mat}
\end{equation}
The gradient of the basis function is computed from the gradients
on the reference cell using (\ref{eq:co_mat})
\[
\nabla\varphi_{L,i}(\vx_{L,j})=\frac{1}{\det\jacob_{L}'(\hat{\vx}_{i})}\text{co}(\jacob_{L}'(\hat{\vx}_{j}))\hat{\nabla}\hat{\varphi}_{i}(\hat{\vx}_{j}).
\]
In the same way, the scaled normal vectors $\v n_{\epsilon}$ on the
faces are computed by the formula
\[
\normal_{\epsilon}(\vx_{L,i})=\text{co}(\jacob_{L}'(\hat{\vx}_{i}))\hat{\normal}_{\epsilon}.
\]
We introduce the following notation for the cofactor matrix\global\long\def\comat{\mathbf{c}}
\[
\comat_{L,i}=\text{co}(\jacob_{L}'(\hat{\vx}_{i})).
\]
The DG scheme then reads
\begin{multline}
\partial_{t}f_{L,i}-\frac{1}{\omega_{L,i}}\sum_{j=0}^{N_{\dorder}-1}\mathbf{v}\cdot\comat_{L,j}\hat{\nabla}\hat{\varphi}_{i}(\hat{\vx}_{j})f_{L,j}\omega_{j}+\\
\frac{1}{\omega_{L,i}}\sum_{\epsilon=1}^{6}\mu_{i}^{\epsilon}\left(\mathbf{v}\cdot\comat_{L,i}\hat{\normal}_{\epsilon}{}^{+}f_{L,i}+\mathbf{v}\cdot\comat_{L,i}\hat{\normal}_{\epsilon}{}^{-}f_{R,i'}\right)=0.\label{eq:DG_reduit}
\end{multline}
On boundary GL points, the value of $f_{R,i'}$ is given by the boundary
condition
\[
f_{R,i'}=f^{b}(\vx_{L,i}),\quad\vx_{L,i}=\vx_{R,i'}.
\]
For practical reasons, it is interesting to also consider $f_{R,i'}$
as an artificial unknown in the fictitious cell. The fictitious unknown
is then a solution of the differential equation
\begin{equation}
\partial_{t}f_{R,i'}=0.\label{eq:fictitious_ODE}
\end{equation}
In the end, if we put all the unknowns in a large vector $\v F(t)$,
(\ref{eq:DG_reduit}), (\ref{eq:fictitious_ODE}) read as a large
system of coupled differential equations\global\long\def\transmat{\mathbf{L}_{h}}
\begin{equation}
\partial_{t}\v F=\transmat\v F.\label{eq:linear_ODE}
\end{equation}

In the following, we call $\transmat$ the transport matrix. The transport
matrix satisfies the following properties:
\begin{itemize}
\item $\transmat\v F=0$ if the components of $\v F$ are all the same.
\item Let $\v F$ be such that the components corresponding to the boundary
term vanish. Then $\v F^{T}\transmat\v F\leq0$. This dissipation
property is a consequence of the choice of an upwind numerical flux
\cite{johnson1984finite}.
\item In many cases, and with a good numbering of the unknowns in $\v F$,
$\transmat$ has a triangular structure. This aspect is discussed
in Subsection \ref{subsec:triangular_subsection}.
\end{itemize}
As stated above, we actually have to apply a transport solver for
each constant velocity $\v v_{i}$.

Let $L$ be a cell of the mesh $\mathcal{M}$ and $\vx_{i}$ a GL
point in $L$. As in the scalar case, we denote by $\vf_{L,i}$ the
approximation of $\vf$ in $L$ at GL point $i$. In the sequel, with
an abuse of notation and according to the context, we may continue
to note $\v F(t)$ the big vector made of all the vectorial values
$\vf_{L,j}$ at all the GL points $j$ in all the (real or fictitious)
cells $L$.

We may also continue to denote by $\v L_{h}$ the matrix made of the
assembly of all the transport operators for all velocities $\v v_{i}$.
With a good numbering of the unknowns it is possible in many cases
to suppose that $\v L_{h}$ is block-triangular. More precisely, because
in the transport step the equations are uncoupled, we see that $\v L_{h}$
can be made block-diagonal, each diagonal block being itself block-triangular.
See Section \ref{subsec:triangular_subsection}.

\subsection{\label{subsec:Palindromic-time-integration}Palindromic time integration}

We can also define an approximation $\v N_{h}$ of the collision operator
$\v N$. We define by $\v F^{eq}(\v F)$ the big vector made of all
the $\vf^{eq}(\vf_{L,i})$, $L\in\mathcal{M},$ $i=0\ldots N_{d}-1.$

We set
\begin{equation}
\v N_{h}\v F=\frac{1}{\tau}(\v F^{eq}(\v F)-\v F).\label{eq:collision_num}
\end{equation}

Similarly we note $\mathbf{G}_{h}$ the discrete approximation of
the kinetic source term $\vg$.

The DG approximation of (\ref{eq:kin_bgk}) finally reads
\[
\partial_{t}\v F=\v L_{h}\v F+\v N_{h}\v F+\mathbf{G}_{h}\mathbf{F}.
\]
We use the following Crank-Nicolson second order time integrator for
the transport equation:
\begin{equation}
\exp(\Delta t\v L_{h})\simeq T_{2}(\Delta t):=(\v I+\frac{\Delta t}{2}\v L_{h})(\v I-\frac{\Delta t}{2}\v L_{h})^{-1}.\label{eq:transport_cn}
\end{equation}
Similarly, for the collision integrator, we use
\[
\exp(\Delta t\v N_{h})\simeq C_{2}(\Delta t):=(\v I+\frac{\Delta t}{2}\v N_{h})(\v I-\frac{\Delta t}{2}\v N_{h})^{-1}.
\]
Because during the collision step, the conservative variables $\vw=\v P\vf$
do not change, the collision integrator is only apparently implicit.
We have the explicit formula:
\begin{equation}
C_{2}(\Delta t)\v F=\frac{(2\tau-\Delta t)\v F}{2\tau+\Delta t}+\frac{2\Delta t\v F^{eq}(\v F)}{2\tau+\Delta t}.\label{eq:collision_cn}
\end{equation}

The source operator is also approximated by a Crank-Nicolson integrator

\[
\exp(\Delta t\v G_{h})\simeq S_{2}(\Delta t):=(\v I+\frac{\Delta t}{2}\v G_{h})(\v I-\frac{\Delta t}{2}\v G_{h})^{-1},
\]

requiring to solve a nonlinear local equation whenever $\vg$ depends
on $\vf$. 

If $\tau>0$, we observe that the operators $T_{2}$ and $C_{2}$
are \emph{time-symmetric}: if we set $O_{2}=T_{2}$ , $O_{2}=C_{2}$,
or $O_{2}=S_{2}$, then $O_{2}$ satisfies
\begin{equation}
O_{2}(-\Delta t)=O_{2}(\Delta t)^{-1},\quad O_{2}(0)=Id.\label{eq:prop_sym}
\end{equation}

This property implies that $O_{2}$ is necessarily a second order
approximation of the exact integrator \cite{mclachlan2002splitting,hairer2006geometric}.
When $\tau=0,$ we also remark that
\[
C_{2}(\Delta t)\v F=2\v F^{eq}(\v F)-\v F\neq\v F
\]
and then $C_{2}$ does not satisfy (\ref{eq:prop_sym}) anymore.

For $\tau>0$, the Strang formula permits us to construct a five steps
second order time-symmetric approximation
\[
M_{2}^{s}(\Delta t)=T_{2}(\frac{\Delta t}{2})S_{2}(\frac{\Delta t}{2})C_{2}(\Delta t)S_{2}(\frac{\Delta t}{2})T_{2}(\frac{\Delta t}{2})=\exp\left(\Delta t\left(\v L_{h}+\v N_{h}+\mathbf{S}_{h}\right)\right)+O(\Delta t^{3}),
\]

and a three step one

\[
M_{2}(\Delta t)=T_{2}(\frac{\Delta t}{2})C_{2}(\Delta t)T_{2}(\frac{\Delta t}{2})=\exp\left(\Delta t\left(\v L_{h}+\v N_{h}\right)\right)+O(\Delta t^{3}),
\]

in the source-less case.

However this formula is no more a second order approximation of (\ref{eq:kin_bgk})
when $\tau\to0$. Indeed, when $\tau=0$
\[
M_{2}(0)\v F=2\v F^{eq}(\v F)-\v F.
\]

As explained in \cite{coulette2016palindromic} it is better to consider
the following method, which remains second order accurate even for
infinitely fast relaxation:
\[
M_{2}^{kin}(\Delta t)=T_{2}(\frac{\Delta t}{4})C_{2}(\frac{\Delta t}{2})T_{2}(\frac{\Delta t}{2})C_{2}(\frac{\Delta t}{2})T_{2}(\frac{\Delta t}{4}).
\]

By palindromic compositions of the second order method $M_{2}^{kin}$
it is then very easy to achieve any even order of accuracy (see \cite{coulette2016palindromic}).
However, in this paper, we concentrate on the parallel optimization
of the method and we shall only present numerical results at second
order for the limit system \ref{eq:macro_limit_system}. To that end
it is sufficient to use the method $M_{2}$, as $PM_{2}(0)$ properly
converges towards identity on the macroscopic variable space when
$\tau\rightarrow0$. 

\section{Optimization of the kinetic solver}

In this section, we describe the optimizations that can be applied
in the implementation of the previous numerical method.

\subsection{\label{subsec:triangular_subsection}Triangular structure of the
transport matrix}

\begin{figure}
\includegraphics[height=6cm]{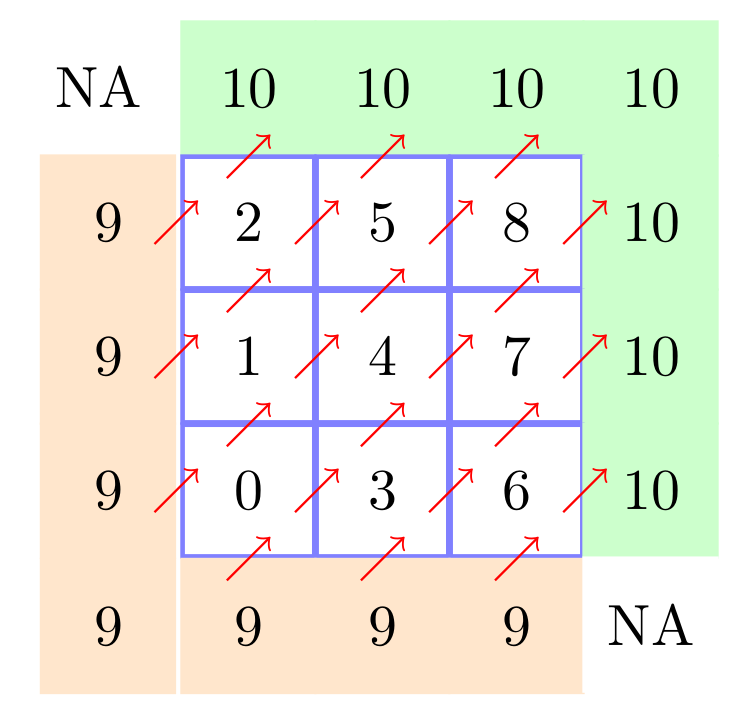}\includegraphics[height=10cm]{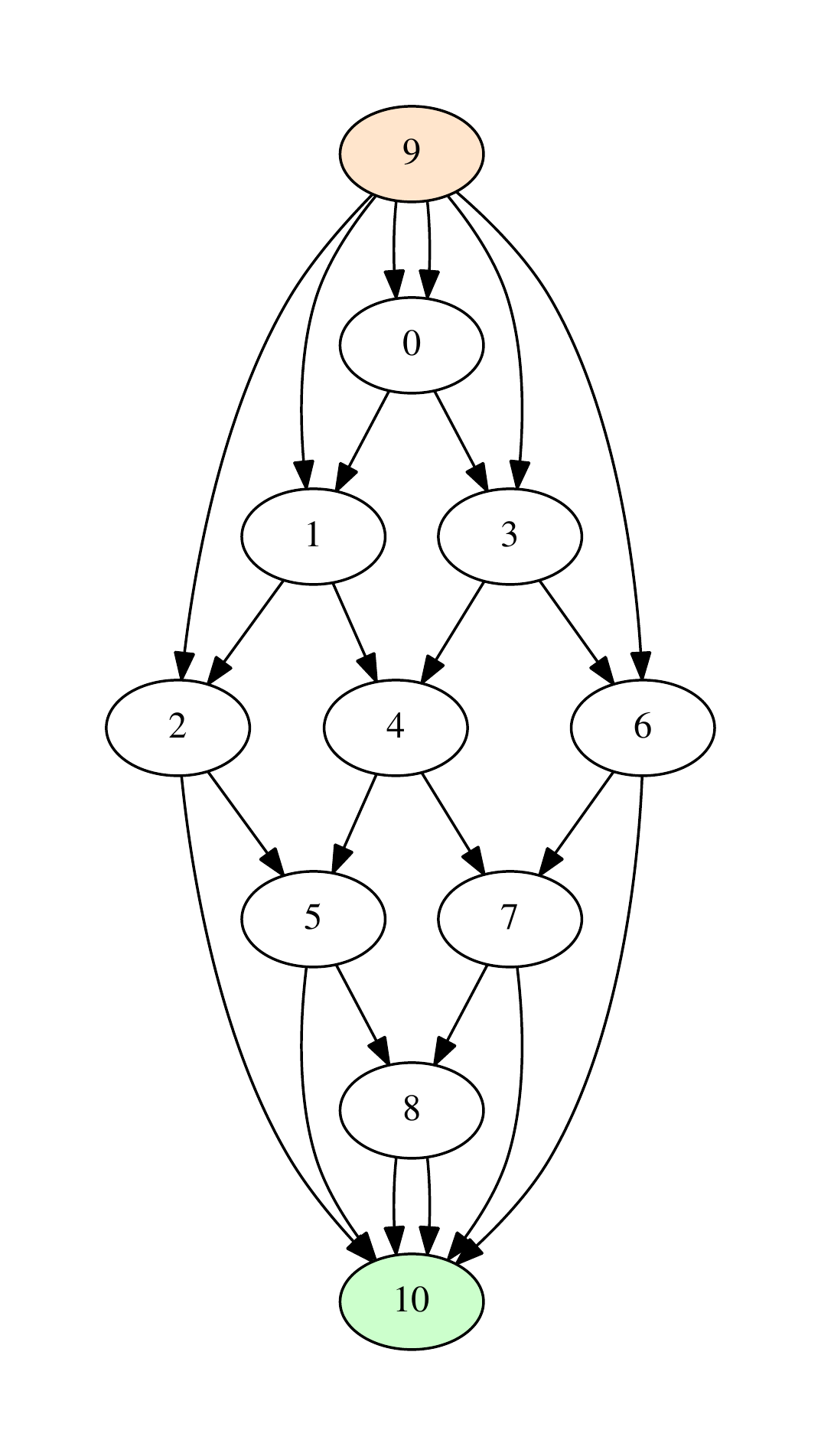}

\caption{\label{fig:dependency_graph}Construction of the dependency graph.
Left: example of mesh (it is structured here but it is not necessary)
with 9 interior cells. The velocity field $v$ is indicated by red
arrows. We add two fictitious cells: one for the upwind boundary condition
(cell $9$) and one for the outflow part of $\partial\Omega$ (cell
$10$). Right: the corresponding dependency graph $\mathcal{G}.$
By examining the dependency graph, we observe that the values of $\protect\v F^{n+1}$
in cell $0$ have to be computed first, using the boundary conditions.
Then cells $1$ and $3$ can be computed in parallel, then cells $2,$
$4,$ and $6$ can be computed in parallel, then \emph{etc}.}
\end{figure}

Because of the upwind structure of the numerical flux, it appears
that the transport matrix is often block-triangular. This is very
interesting because this allows to apply implicit schemes to (\ref{eq:linear_ODE})
without the costly inversion of linear systems \cite{moustafa20153d}.
We can provide the formal structure of $\transmat$ through the construction
of a directed graph $\mathcal{G}$ with a set of vertices $\mathcal{V}$
and a set of edges $\mathcal{E}\subset\mathcal{V}\times\mathcal{V}$.
The vertices of the graph are associated to the (real or fictitious)
cells of $\mathcal{M}$. Consider now two cells $L$ and $R$ with
a common face $F_{LR}$. We denote by $\normal_{LR}$ the normal vector
on $F_{LR}$ oriented from $L$ to $R$. If there is at least one
GL point $\vx$ on $F_{LR}$ such that 
\[
\normal_{LR}(\vx)\cdot\mathbf{v}>0,
\]
then the edge from $L$ to $R$ belongs to the graph:
\[
(L,R)\in\mathcal{E},
\]
see Figure \ref{fig:dependency_graph}.

In (\ref{eq:DG_reduit}) we can distinguish between several kinds
of terms. We write
\[
\partial_{t}f_{L}+\Gamma_{L\leftarrow L}f_{L}+\sum_{(R,L)\in\mathcal{E}}\Gamma_{L\leftarrow R}f_{R},
\]
with
\begin{multline*}
\Gamma_{L\leftarrow L}f_{L}=-\frac{1}{\omega_{L,i}}\sum_{j=0}^{N_{\dorder}-1}\mathbf{v}\cdot\comat_{L,j}\hat{\nabla}\hat{\varphi}_{i}(\hat{\vx}_{j})f_{L,j}\omega_{j}+\\
\frac{1}{\omega_{L,i}}\sum_{\epsilon=1}^{6}\mu_{i}^{\epsilon}\mathbf{v}\cdot\comat_{L,i}\hat{\normal}_{\epsilon}{}^{+}f_{L,i},
\end{multline*}
and, if $(R,L)\in\mathcal{E}$,
\[
\Gamma_{L\leftarrow R}f_{R}=\frac{1}{\omega_{L,i}}\mu_{i}^{\epsilon}\mathbf{v}\cdot\comat_{L,i}\hat{\normal}_{\epsilon}{}^{-}f_{R,i'}.
\]
We can use the following convention
\begin{equation}
(R,L)\notin\mathcal{E}\Rightarrow\Gamma_{L\leftarrow R}=0.\label{eq:convention_downwind}
\end{equation}
$\Gamma_{L\leftarrow L}$ contains the terms that couple the values
of $f$ inside the cell $L$. They correspond to diagonal blocks of
size $(\dorder+1)^{\ddim}\times(\dorder+1)^{\ddim}$ in the transport
matrix $\transmat$. $\Gamma_{L\leftarrow R}$ contains the terms
that couple the values inside cell $L$ with the values in the neighboring
upwind cell $R$. If $R$ is a downwind cell relatively to $L$ then
$\mu_{i}^{\epsilon}\mathbf{v}\cdot C_{L,i}\hat{\normal}_{\epsilon}{}^{-}=0$
and $\Gamma_{L\leftarrow R}=0$ is indeed compatible with the above
convention (\ref{eq:convention_downwind}).

Once the graph $\mathcal{G}$ is constructed, we can analyze it with
standard tools. If it contains no cycle, then it is called a Directed
Acyclic Graph (DAG). Any DAG admits a topological ordering of its
nodes. A topological ordering is a numbering of the cells $i\mapsto L_{i}$
such that if there is a path from $L_{i}$ to $L_{j}$ in $\mathcal{G}$
then $j>i$. In practice, it is useful to remove the fictitious cells
from the topological ordering. In our implementation they are put
at the end of the list.

Once the new ordering of the graph vertices is constructed, we can
construct a numbering of the components of $\v F$ by first numbering
the unknowns in $L_{0}$ then the unknowns in $L_{1}$, \emph{etc}.
More precisely, we set
\[
F_{kN_{\dorder}+i}=f_{L_{k},i}.
\]
Then, with this ordering, the matrix $\transmat$ is lower block-triangular
with diagonal blocks of size $(\dorder+1)^{\ddim}\times(\dorder+1)^{\ddim}$.
It means that we can apply implicit schemes to (\ref{eq:linear_ODE})
without inverting large linear systems.

As stated above, we actually have to apply a transport solver for
each constant velocity $\v v_{i}$. In the sequel, with another abuse
of notation and according to the context, we continue to note $\v F$
the big vector made of all the vectorial values $\vf_{L,j}$ at all
the GL points $j$ in all the (real or fictitious) cells $L$.

We may also continue to denote by $\v L_{h}$ the matrix made of the
assembly of all the transport operators for all velocities $\v v_{i}$.
With a good numbering of the unknown it is still possible to suppose
that $\v L_{h}$ is block-triangular. More precisely, as in the transport
step the equations are uncoupled, we see that $\v L_{h}$ can be made
a block-diagonal matrix, each diagonal block being itself block-triangular.

\subsection{Parallelization of the implicit solver}

In this section, we explain how it is possible to parallelize the
transport solver. Here again we consider the single transport equation
(\ref{eq:single_transport}) and the associated differential equation
(\ref{eq:linear_ODE}). We apply a second order Crank-Nicolson implicit
scheme. We have explained in Section \ref{subsec:Palindromic-time-integration}
how to increase the order of the scheme. We compute an approximation
$\v F^{n}$ of $\v F(n\Delta t$). The implicit scheme reads
\begin{equation}
(\v I-\Delta t\transmat)\v F^{n+1}=(\v I+\Delta t\transmat)\v F^{n}.\label{eq:crank_nicolson}
\end{equation}
As explained above, the matrices $(\v I-\Delta t\transmat)$ and $(\v I+\Delta t\transmat)$
are lower triangular. It is thus possible to solve the linear system
explicitly cell after cell, assuming that the cells are numbered in
a topological order.

It is possible to perform further optimization by harnessing the parallelism
exhibited by the dependency graph. Indeed, once the values of $f$
in the first cell are computed, it is generally possible to compute
in parallel the values of $f$ in neighboring downwind cells. For
example, as can be seen on Figure \ref{fig:dependency_graph}, once
the values in cells 0, 1 and 2 are known, we can compute independently,
and in parallel, the values in cells 2, 4 and 6.

We observe that at the beginning and at the end of the time step,
the computations are ``less parallel'' than in the middle of the
time step, where the parallelism is maximal.

Implementing this algorithm with OpenMP or using pthread is not very
difficult. However, it requires to compute the data dependencies between
the computational tasks carefully, and to set adequate synchronization
points in order to get correct results. In addition, a rough implementation
will probably not exhibit optimized memory access. Therefore, we have
decided to rely on a more sophisticated tool called StarPU\footnote{\url{http://starpu.gforge.inria.fr}}
for submitting the parallel tasks to the available computational resources. 

StarPU is a runtime system library developed at Inria Bordeaux \cite{AugAumFurNamThi2012EuroMPI}.
It relies on the data-based parallelism paradigm.

The user has first to split its whole problem into elementary computational
tasks. The elementary tasks are then implemented into \textit{codelets},
which are simple C functions. The same task can be implemented differently
into several codelets. This allows the user to harness special acceleration
devices, such as vectorial CPU cores, GPUs or Intel KNL devices, for
example. In the StarPU terminology these devices are called \textit{workers}.

For each task, the user has also to describe precisely what are the
input data, in \textit{read} mode, and the output data, in \textit{write}
or \textit{read-write} mode. The user then submits the task in a sequential
way to the StarPU system. StarPU is able to construct at runtime a
task graph from the data dependencies. The task graph is analyzed
and the tasks are scheduled automatically to the available workers
(CPU cores, GPUs, \emph{etc}.). If possible, they are executed in
parallel into concurrent threads. The data transfer tasks between
the threads are automatically generated and managed by StarPU, which
greatly simplifies the programming.

When a StarPU program is executed, it is possible to choose among
several schedulers. The simplest \textit{eager} scheduler adopts a
very simple strategy, where the tasks are executed in the order of
submission by the free workers, without optimization. More sophisticated
schedulers, such as the \textit{dmda} scheduler, are able to measure
the efficiency of the different codelets and the data transfer times,
in order to apply a more efficient allocation of tasks.

Recently a new data access mode has been added to StarPU: the \textit{commute}
mode. In a task, a buffer of data can now be accessed in commute mode,
in addition to the write or read-write modes. A commute access tells
to StarPU that the execution of the corresponding task may be executed
before or after other tasks containing commutative access. This allows
StarPU to perform additional optimizations.

There exists also a MPI version of StarPU. In the MPI version, the
user has to decide an initial distribution of data among the MPI nodes.
Then the tasks are submitted as usual (using the function starpu\_mpi\_insert\_task
instead of starpu\_insert\_task). Required MPI communications are
automatically generated by StarPU. For the moment, this approach does
not guarantee a good load balancing. It is the responsibility of the
user to migrate data from one MPI node to another for improving the
load balancing, if necessary.

\subsection{\label{subsec:Macrocell-approach}Macrocell approach}

StarPU is quite efficient, but there is an unavoidable overhead due
to the task submissions and to the on-the-fly construction and analysis
of the task graph. Therefore it is important to ensure that the computational
tasks are not too small, in which case the overhead is not amortized,
or not too big, in which case some workers are idle.

For achieving the right balance, we have decided not to apply directly
the above task submission algorithm to the cells but to groups of
cells instead.

The implementation of the whole kinetic scheme has been made into
the \texttt{schnaps} software\footnote{\selectlanguage{french}%
\url{http://schnaps.gforge.inria.fr/}\selectlanguage{english}%
}. \texttt{schnaps} is a C99 software dedicated to the numerical simulation
of conservation laws. 

\begin{figure}
\begin{centering}
\includegraphics[height=8cm]{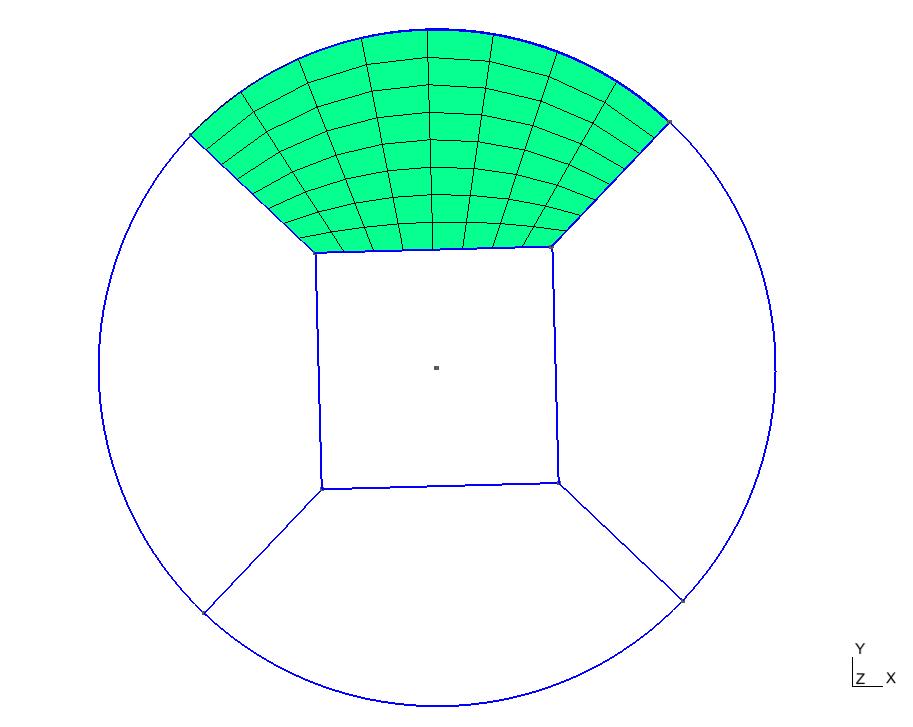}
\par\end{centering}
\caption{\label{fig:macrocell_mesh}Macrocell approach: an example of a mesh
made of five macrocells. Each macrocell is then split into several
subcells. Only the subcells of the top macrocell are represented here
(in green).}
\end{figure}

In \texttt{schnaps} we construct first a \textit{macromesh} of the
computational domain. Then each \textit{macrocell} of the macromesh
is split into subcells. See Figure \ref{fig:macrocell_mesh}. We also
arrange the subcells into a regular sub-mesh of the macrocells. In
this way, it is possible to apply additional optimizations. For instance,
the subcells $L$ of a same macrocell $\mathcal{L}$ can now share
the same geometrical transformation $\jacob_{L}$, which saves memory.

In \texttt{schnaps} we have also defined an \textit{interface} structure
in order to manage data communications between the macrocells. An
interface contains the faces that are common to two neighboring macrocells.
We do not proceed exactly as in Section \ref{subsec:triangular_subsection}
where the vertices of graph $\mathcal{G}$ were associated to cells
and the edges to faces. Instead, we construct an upwind graph whose
vertices are associated to macrocells, and edges to interfaces. This
graph is then sorted, and the macrocells are numbered in a topological
order.

For solving one time step of one transport equation (\ref{eq:crank_nicolson}),
we split the computations into several elementary operations: for
all macrocell $\mathcal{L}$ taken in a topological order, we perform
the following tasks:
\begin{enumerate}
\item Volume residual assembly: this task computes in the macrocell $\mathcal{L}$
the part of the right hand side of (\ref{eq:crank_nicolson}) that
comes from the values of $f$ inside $\mathcal{L}$;
\item Interface residual assembly: this task computes, in the macrocell
$\mathcal{L}$, the part of the right hand side of (\ref{eq:crank_nicolson})
that comes from upwind interface values;
\item Boundary residual assembly: this task computes, in the macrocell $\mathcal{L}$,
the part of the right hand side of (\ref{eq:crank_nicolson}) that
comes from upwind boundaries values.
\item \label{enu:vol-solve}Volume solve: this task solves the local transport
linear system in the macrocell.
\item Extraction: this task copies the boundary data of $\mathcal{L}$ to
the neighbor downwind interfaces.
\end{enumerate}
Let us point out that in step \ref{enu:vol-solve} above, the macrocell
local transport solver is reassembled and refactorized at each time
step: we have decided not to store a sparse linear system in the macrocell
for each velocity $\v v_{i}$, in order to save memory. The local
sparse linear system is solved thanks to the KLU library \cite{davis2010algorithm}.
This library is able to detect efficiently sparse triangular matrix
structures, which makes the resolution quite efficient. In practice,
the factorization and resolution time of the KLU solver is of the
same order as the residual assembly time.

In schnaps, we use the MPI version of StarPU. The macromesh is initially
split into several subdomains and the subdomains are distributed to
the MPI nodes. Then the above tasks are launched asynchronously with
the starpu\_mpi\_insert\_task function. MPI communications are managed
automatically by StarPU.

It is clear that if we were solving a single transport equation our
strategy would be very inefficient. Indeed, the downwind subdomains
would have to wait for the end of the computations of the upwind subdomains.
We are saved by the fact that we have to solve many transport equations
in different directions. This helps the MPI nodes to be equally occupied.
Our approach is more efficient if we avoid a domain decomposition
with internal subdomains, because these subdomains have to wait the
results coming from the boundaries.

In our approach it is also essential to launch the tasks in a completely
asynchronous fashion. In this way, if a MPI node is waiting for results
of other subdomains for a given velocity $\v v_{i}$ it is not prevented
from starting the computation for another velocity $\v v_{j}.$

\subsection{Collisions}

In this section we explain how is computed the collision step (\ref{eq:collision_cn}).
The computations are purely local to each GL point, which makes the
collision step embarrassingly parallel. However it is not so obvious
to attain high efficiency because of memory access. If the values
of $\v F$ are well arranged in memory in the transport stage, it
means that the values of $\vf$ attached to a given velocity $\v v_{i}$
are close in memory, for a better data locality. On the contrary,
in the collision step at a given GL point, a better locality is achieved
if the values of $\vf$ corresponding to different velocities are
close in memory. Additional investigations and tests are needed in
order to evaluate the importance of data locality in our algorithm.

In our implementation, we have adopted the following strategy. We
have first identified the following task:
\begin{enumerate}
\item Reduction task for a velocity $\v v_{i}$: this task is associated
with one macrocell $\mathcal{L}$. It computes the contribution to
$\vw$ of the components of $\v f$ that have been transported at
velocity $\v v_{i}$ with formula (\ref{eq:micro_to_macro}). The
StarPU access to the buffer containing $\vw$ is performed in read-write
and commute modes. In this way the contribution from each velocity
can be added to $\vw$ as soon as it is available.
\item \label{enu:relax_step_iv}Relaxation task for a velocity $\v v_{i}$:
this task is associated to one macrocell $\mathcal{L}$. Once $\vw$
is known, it computes the components of $\vf^{eq}$ corresponding
to velocity $\v v_{i}.$ Then it computes the relaxation step (\ref{eq:collision_cn})
for the associated component of $\v f$.
\end{enumerate}
In step \ref{enu:relax_step_iv} we can separate the computations
for each velocity because the collision term (\ref{eq:collision_num})
is diagonal. Some Lattice Boltzmann Methods rely on non-diagonal relaxations.
It can be useful for representing more general viscous terms for instance.
For non-diagonal relaxation we would have to change a little the algorithm.

We can now make a few comments about the storage cost of the method.
In the end, we have to store at each GL point $\v x_{i}$ and each
cell $L$ the values of $\vf_{L,i}$ and $\vw_{L,i}$. We do not have
to keep the values of the previous time-step, $\vf_{L,i}$ and $\vw_{L,i}$
can be replaced by the new values as soon as they are available. In
this sense, our method is ``low-storage''. As explained in \cite{coulette2016palindromic}
it is also possible to increase the time order of the method, without
increasing the storage.

\subsection{Scaling test}

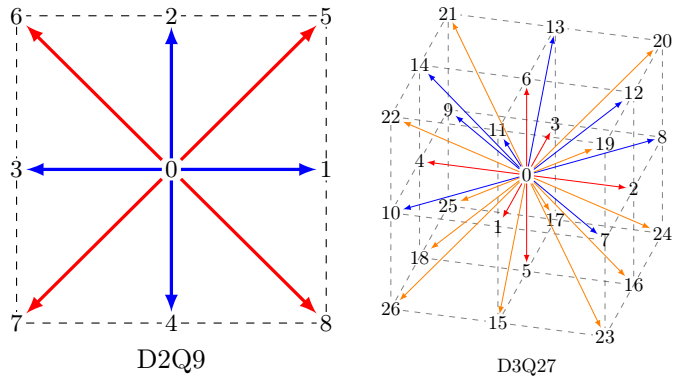
\begin{figure}
\resizebox{!}{5cm}{
\begin{tikzpicture}[scale=1.0]
\begin{small}
\def\vscale{2.0};

\coordinate (C0) at (0,0);
\coordinate (C1) at (\vscale,0);
\coordinate (C2) at (0,\vscale);
\coordinate (C3) at (-\vscale,0);
\coordinate (C4) at (0,-\vscale);
\coordinate (C5) at (\vscale,\vscale);
\coordinate (C6) at (-\vscale,\vscale);
\coordinate (C7) at (-\vscale,-\vscale);
\coordinate (C8) at (\vscale,-\vscale);
\foreach \p in {0,1,...,8} {
	\draw (C\p) node[inner sep=1pt] (D\p){\p};
};
\draw[dashed] (D8)--(D1);
\draw[dashed] (D1)--(D5);
\draw[dashed] (D5)--(D2);
\draw[dashed] (D2)--(D6);
\draw[dashed] (D6)--(D3);
\draw[dashed] (D3)--(D7);
\draw[dashed] (D7)--(D4);
\draw[dashed] (D4)--(D8);

\foreach \c [count = \p] in {blue,blue,blue,blue,red,red,red,red}{
	\draw[-latex,very thick,\c] (D0)--(D\p);
};
\end{small}
\coordinate (T) at (0,-2.5);
\draw (T) node{D2Q9};

\end{tikzpicture}
}\resizebox{!}{5cm}{
\tdplotsetmaincoords{61}{105}
\begin{tikzpicture}[tdplot_main_coords,scale=1.0]
\def\v{2.0};
\coordinate (C0) at (0,0,0);
\coordinate (C1) at (\v,0,0);
\coordinate (C2) at (0,\v,0);
\coordinate (C3) at (-\v,0,0);
\coordinate (C4) at (0,-\v,0);
\coordinate (C5) at (0,0,-\v);
\coordinate (C6) at (0,0,\v);

\coordinate (C7) at (\v,\v,0);
\coordinate (C8) at (-\v,\v,0);
\coordinate (C9) at (-\v,-\v,0);
\coordinate (C10) at (\v,-\v,0);
\coordinate (C11) at (\v,0,\v);
\coordinate (C12) at (0,\v,\v);
\coordinate (C13) at (-\v,0,\v);
\coordinate (C14) at (0,-\v,\v);
\coordinate (C15) at (\v,0,-\v);
\coordinate (C16) at (0,\v,-\v);
\coordinate (C17) at (-\v,0,-\v);
\coordinate (C18) at (0,-\v,-\v);

\coordinate (C19) at (\v,\v,\v);
\coordinate (C20) at (-\v,\v,\v);
\coordinate (C21) at (-\v,-\v,\v);
\coordinate (C22) at (\v,-\v,\v);
\coordinate (C23) at (\v,\v,-\v);
\coordinate (C24) at (-\v,\v,-\v);
\coordinate (C25) at (-\v,-\v,-\v);
\coordinate (C26) at (\v,-\v,-\v);
\foreach \p in {0,1,...,26} {
	\draw (C\p) node[inner sep=1pt] (D\p){\p};
};

\draw[dashed,opacity=0.5] (D19)--(D12);
\draw[dashed,opacity=0.5] (D12)--(D20);

\draw[dashed,opacity=0.5] (D20)--(D13);
\draw[dashed,opacity=0.5] (D13)--(D21);
\draw[dashed,opacity=0.5] (D21)--(D14);
\draw[dashed,opacity=0.5] (D14)--(D22);

\draw[dashed,opacity=0.5] (D22)--(D19);

\draw[dashed,opacity=0.5] (D23)--(D24);
\draw[dashed,opacity=0.5] (D24)--(D25);
\draw[dashed,opacity=0.5] (D25)--(D26);
\draw[dashed,opacity=0.5] (D26)--(D15);
\draw[dashed,opacity=0.5] (D15)--(D23);

\draw[dashed,opacity=0.5] (D7)--(D8);
\draw[dashed,opacity=0.5] (D8)--(D9);
\draw[dashed,opacity=0.5] (D9)--(D4);
\draw[dashed,opacity=0.5] (D4)--(D10);
\draw[dashed,opacity=0.5] (D10)--(D7);
\draw[dashed,opacity=0.5] (D12)--(D6);
\draw[dashed,opacity=0.5] (D6)--(D14);

\draw[dashed,opacity=0.5] (D14)--(D4);
\draw[dashed,opacity=0.5] (D4)--(D18);

\draw[dashed,opacity=0.5] (D18)--(D5);
\draw[dashed,opacity=0.5] (D5)--(D16);

\draw[dashed,opacity=0.5] (D16)--(D12);
\draw[dashed,opacity=0.5] (D19)--(D23);
\draw[dashed,opacity=0.5] (D20)--(D8);
\draw[dashed,opacity=0.5] (D8)--(D24);
\draw[dashed,opacity=0.5] (D21)--(D25);
\draw[dashed,opacity=0.5] (D22)--(D10);
\draw[dashed,opacity=0.5] (D10)--(D26);

\draw[dashed,opacity=0.5] (D11)--(D15);
\draw[dashed,opacity=0.5] (D13)--(D17);
\draw[dashed,opacity=0.5] (D11)--(D6);
\draw[dashed,opacity=0.5] (D6)--(D13);

\draw[dashed,opacity=0.5] (D15)--(D5);
\draw[dashed,opacity=0.5] (D5)--(D17);

\def\mycolora{red}
\def\mycolorb{blue}
\def\mycolorc{orange}
\def\arrowcolors{\mycolora,\mycolora,\mycolora,\mycolora,\mycolora,\mycolora,\mycolorb,\mycolorb,\mycolorb,\mycolorb,\mycolorb,\mycolorb,\mycolorb,\mycolorb,\mycolorc,\mycolorc,\mycolorc,\mycolorc,\mycolorc,\mycolorc,\mycolorc,\mycolorc,\mycolorc,\mycolorc,\mycolorc,\mycolorc}
\foreach \c [count = \i] in \arrowcolors{
	pgfmatheval{\i+1}
	\draw[-latex,\c] (D0)--(D\pgfmathresult);
};
%
\coordinate (T) at (0,0,-4);
\draw (T) node{D3Q27};
\end{tikzpicture}
} 

\caption{\label{fig:D2Q9_D3Q27_nodes}D2Q9 and D3Q27 velocity grids.}
\end{figure}

For all performance tests presented in this section, we used standard
models from the family Lattice-Boltzmann-Method (LBM) kinetic models
devised for the simulation of Euler/Navier Stokes systems. We will
not give a detailed description of their properties from the modeling
point of view: we simply take them as good representative of the typical
workload of kinetic relaxation schemes. The most relevant feature
impacting the performance of our algorithm is the discrete velocity
set of the kinetic model, which determines the task graph structure
of the transport step when combined with a particular mesh topology.
In standard LBM models, velocity sets are usually built-up from a
sequence of pairs of opposite velocities with an additional zero velocity
node. On Figure \ref{fig:D2Q9_D3Q27_nodes} we show the two representative
velocity sets of the $2DQ9$ and $3DQ27$ LBM models. 

\subsubsection{Multithread performance (D2Q9, D3Q{*})}

\begin{figure}
\begin{centering}
\includegraphics[height=8cm]{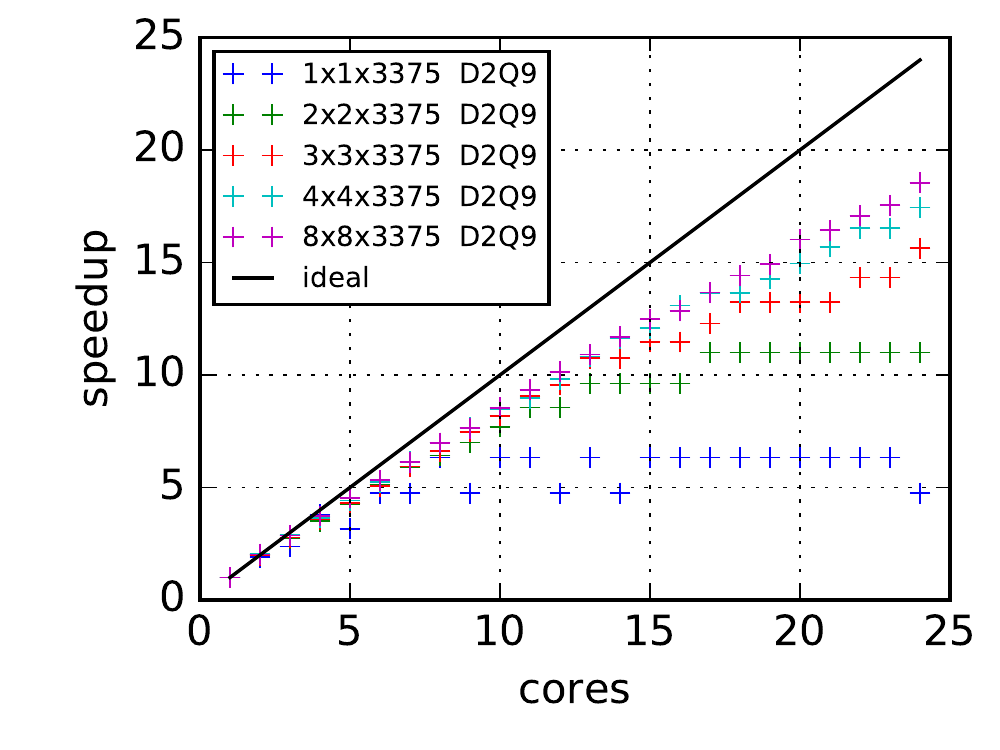}
\par\end{centering}
\caption{\label{fig:D2Q9_speedup_multithread}Multithread scaling for the D2Q9
model and a collection of square meshes from 1 to 64 macro-elements. }
\end{figure}

\begin{figure}
\begin{centering}
\includegraphics[height=8cm]{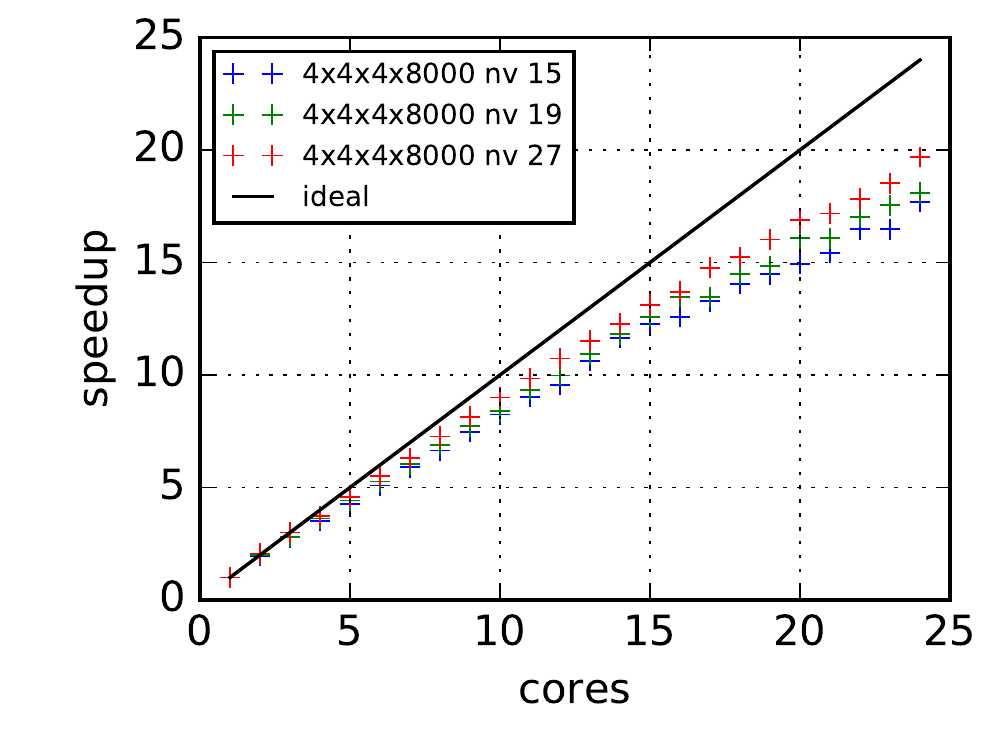}
\par\end{centering}
\caption{\label{fig:D3Qq_speedup_multithread}Multithread scaling on a 4x4x4
macro-element mesh for models D3Q15, D3Q19 and D3Q27}
\end{figure}

We first tested the multithread performance of our implementation
for the full (transport + relaxation) scheme for the standard D2Q9
model. All tests were performed on a single node of the IRMA-ATLAS
cluster, with 24 available cores. We considered several square meshes
build-up from 1 to 64 macrocells. The number of geometric degrees
of freedom per element has been kept constant with a value of $3375$
points per macrocell, so that the workload per macrocell did not change.
For each mesh, we allowed StarPu to use from 1 to the full 24 cores
of the node and measured the total wall time. The results for this
first batch of performance measurements are given in Fig. \ref{fig:D2Q9_speedup_multithread}.
First we verify that for 1 unique macrocell, parallel performance
saturates when the number of cores roughly equals the number of velocities
in the model. This is to be expected, as no topological parallelism
can be exploited in that case. Increasing the number of macro-element
allows to take advantage of topological parallelism. For that workload,
parallel efficiency saturates at about $80
$, which is quite good. On Fig. \ref{fig:D3Qq_speedup_multithread},
we considered on the same cubic mesh three different models differing
by the number of velocity values. Those cases exhibit a large amount
of potential parallelism, due to the large number of velocities combined
with the macrocell decomposition. On an ideal machine, they could
in theory scale perfectly up to 24 cores. The observed saturation,
still around $80
$ efficiency, is still quite good and comes from the unavoidable concurrency
in memory access between the various cores and the scheduling overhead.
It is important to note when considering those results that the bulk
of the computational cost occurs in the transport step of the algorithm:
though the collision step forces synchronization between all the fields
corresponding to individual velocities for the computation of the
macroscopic fields, its actual cost is negligible in regard of the
transport step. 

\subsubsection{MPI scaling: D3Q15 in a torus}

Having verified the good multithread performance of our code on a
single node, we now check whether for larger problem sizes the workload
can be distributed among several computing nodes. To that end, accounting
for the fact that we aim notably at performing simulations for Tokamak
physics, we considered a toroidal mesh subdivided into $720$ macrocells.
The workload distribution across nodes was made using a standard domain
decomposition approach: the mesh was partitioned statically into as
many sub-domains as computing nodes, ranging from $1$to $4$ for
our experiment on the IRMA-ATLAS cluster. From an implementation point
of view, the transition from a multithreaded code to a hybrid MPI/multithread
one is made fairly easy by StarPU. When declaring data to StarPU,
one simply has to specify the MPI process owning the data. At runtime,
each MPI process hosts a local scheduler instance which acts only
on data relevant to the local execution graph. All MPI communications
are handled transparently by the local scheduler when inter-node data
transfers are necessary. In table \ref{tab:mpi_scaling} we show the
wall time for a hundred iterations of the full scheme for the D3Q15
model. The number of available threads per node was set to $14,$
matching the number of velocities actually participating in transport
(there is one null velocity in the model). We observe a super-linear
scaling when the load is spread from $1$to $4$ node. This is not
surprising for such an experiment with fixed total problem size as
both the memory load and size of the local task graph for each decrease
when the number of sub-domains increases.

\begin{figure}
\begin{centering}
\includegraphics[width=4cm]{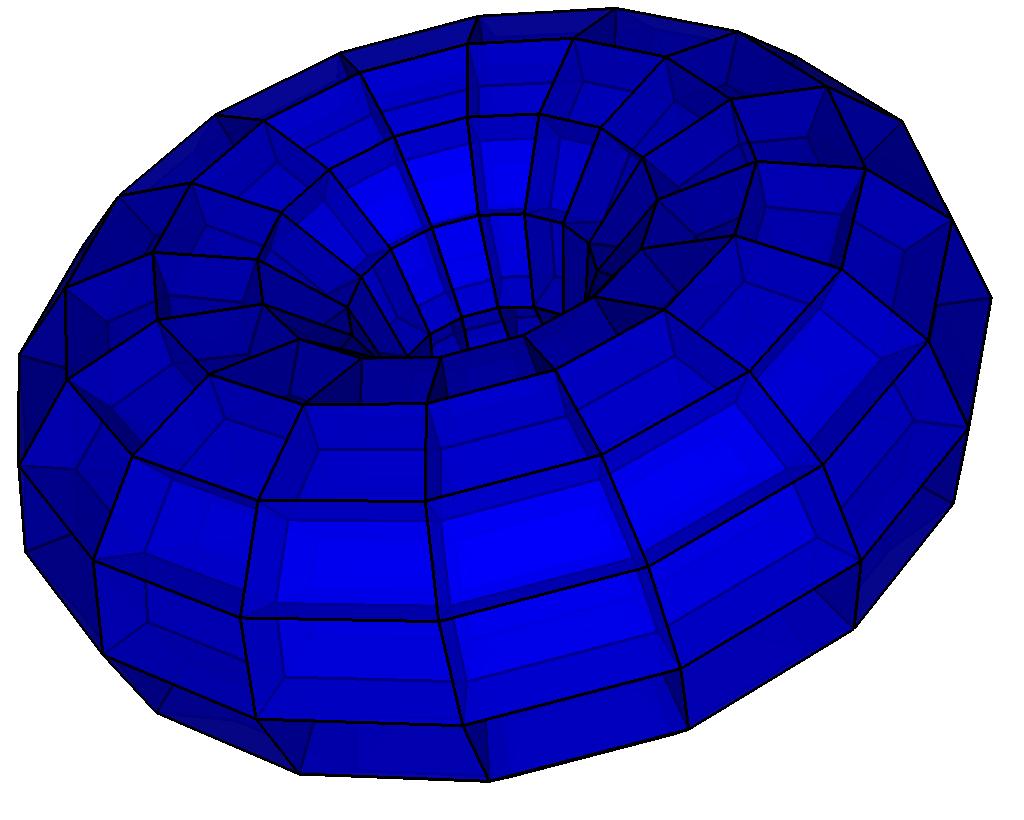}\includegraphics[width=4cm]{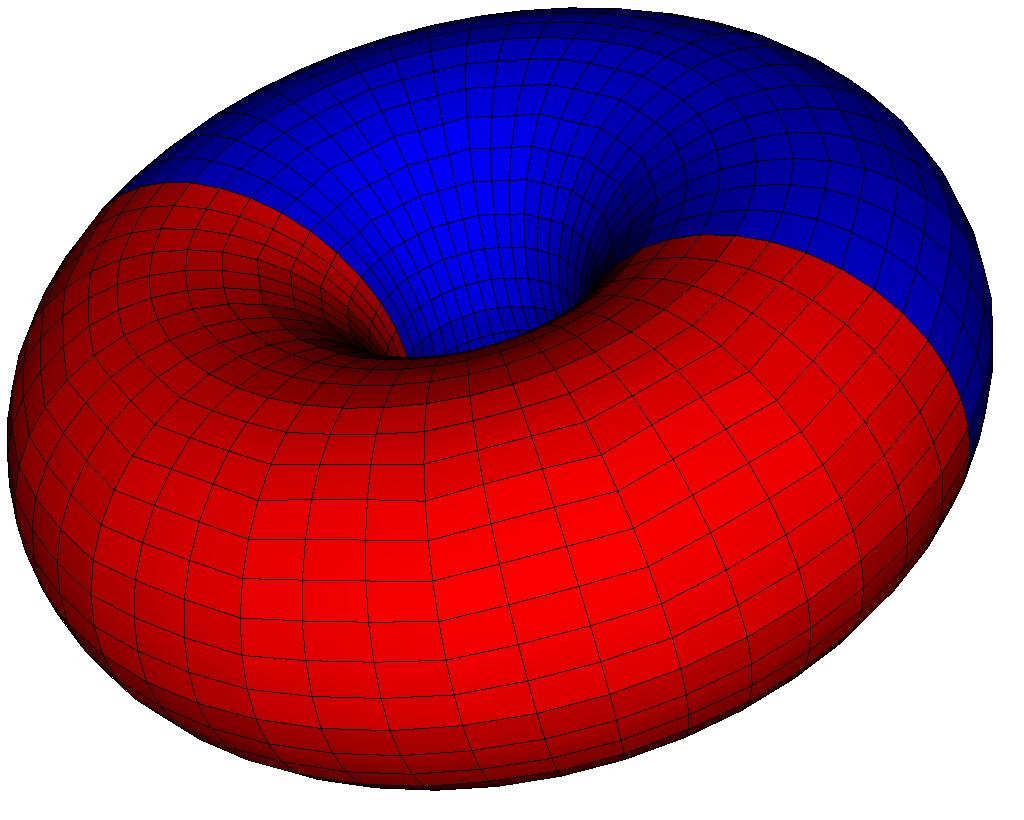}\includegraphics[width=4cm]{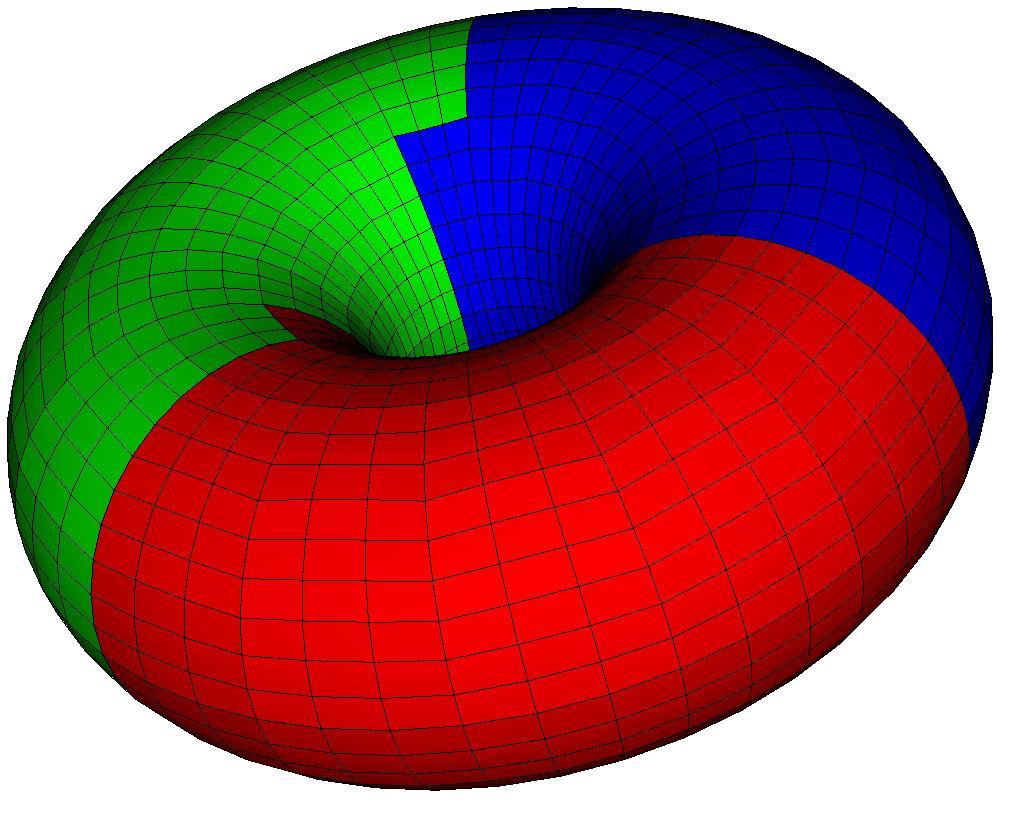}\includegraphics[width=4cm]{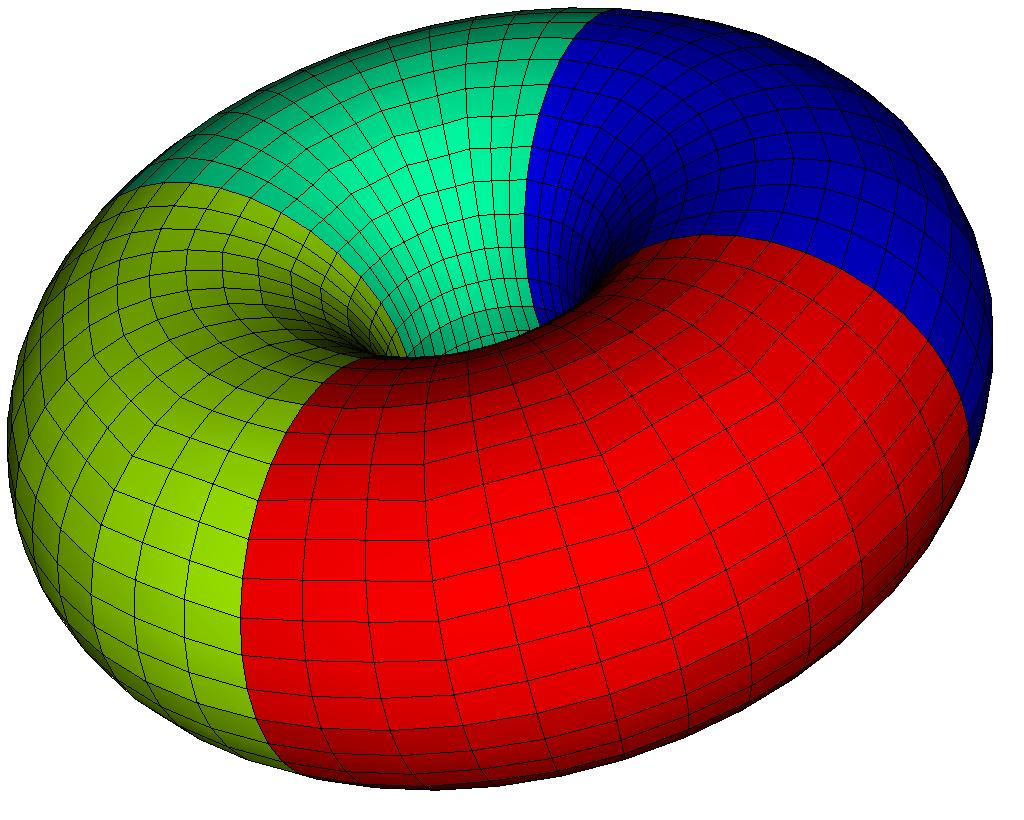}
\par\end{centering}
\caption{\label{fig:torus_partition}Toroidal macromesh ($720$ macrocells)
- Mesh partitions used in the MPI scaling tests. }
\end{figure}

\begin{table}
\begin{center}
\begin{tabular}{|c|cccc|} \hline Nthreads/Nmpi &  1 & 2 & 3 & 4 \\ \hline 14 & 6862 & 2772 &  1491 & 1014  \\ \hline \end{tabular} 
\end{center}

\caption{\label{tab:mpi_scaling}Wall time (in seconds) for the D3Q15 model
for 1 to 4 mpi processes with 14 threads per process.}
\end{table}

\section{Numerical results}

\subsection{Euler with gravity}

For this test case we consider the isothermal Euler equations in two
dimensions with a constant gravity source term

\begin{equation}
\partial_{t}\rho+\partial_{k}(\rho\mathbf{\Vu^{k})}=0,
\end{equation}

\begin{equation}
\partial_{t}(\rho\Vu^{k})+\partial_{j}\boldsymbol{\Pi^{kj}}=\rho\mathbf{g},
\end{equation}

with $\boldsymbol{\Pi}=\left[\begin{array}{cc}
\rho c^{2}+\rho u_{x}^{2} & \rho u_{x}u_{y}\\
\rho u_{x}u_{y} & \rho c^{2}+u_{y}^{2}
\end{array}\right]$ and $\mathbf{g}=-g\mathbf{e}_{y}$. 

The conservative variables vector is thus $\vw=[\rho,\rho u_{x},\rho u_{y}]^{t}$
. The kinetic model is the standard $D2Q9$ one with nine velocities 

\begin{equation}
\vV=\lambda\mathrm{\textrm{diag}}\left[(0,0),(1,0),(0,1),(-1,0),(0,-1),(1,1),(-1,1),(-1,-1),(1,-1)\right]
\end{equation}

and the $(3\times9)$ projection matrix $P$ reads

\begin{equation}
P=\begin{bmatrix}1 & 1 & 1 & 1 & 1 & 1 & 1 & 1 & 1\\
0 & \lambda & 0 & -\lambda & 0 & \lambda & -\lambda & -\lambda & \lambda\\
0 & 0 & \lambda & 0 & -\lambda & \lambda & \lambda & -\lambda & -\lambda
\end{bmatrix},
\end{equation}

i.e $\rho=\sum_{i}f_{i}$, $\rho\Vu=\sum_{i}f_{i}\Vvi$. 

The equilibrium distribution function is given by 

\begin{equation}
f_{i}=w_{i}\rho\left(1+\frac{(\Vu\cdot\Vvi)}{c^{2}}+\frac{(\Vu\cdot\Vvi)^{2}}{2c^{4}}-\frac{\Vu\cdot\Vu}{2c^{2}}\right)
\end{equation}

with $c=\lambda/\sqrt{3}$, and the weights $w_{0}=\frac{4}{9}$,
$w_{i}=\frac{1}{9}$ for $i=1,\dots,4$, $w_{i}=\frac{1}{36}$ for
$i=5,\dots,8$.

The stationary solution for a fluid at rest in the gravity field is

\begin{equation}
\rho=\rho_{0}\exp(-gy/c^{2}),\qquad\Vu=0
\end{equation}

For this test case, the numerical scheme is made up of three stages:
a transport step (T), a source step (S) where the source is applied
on the equilibrium part of the distribution function, and the collision
step (C). Due to the absence of explicit time dependency and the linearity
in $\vw$ of the source, the local nonlinear resolution of the source
operator converges in one Picard iteration. All steps are implemented
as weighted implicit schemes, parametrized by a weight $\theta$ and
a time step$\Delta t$.  

\begin{figure}
\begin{centering}
\includegraphics[height=8cm]{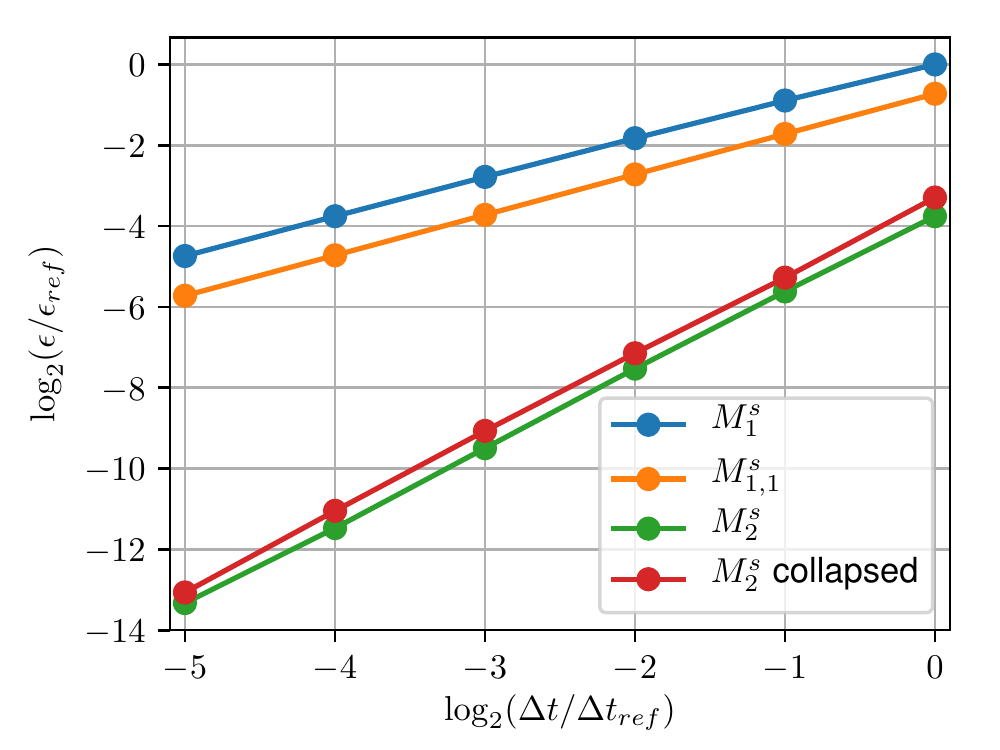}
\par\end{centering}
\caption{\label{fig:euler_gravity_order}Time order convergence for the $2D$
Euler gravity test case with $D2Q9$ model. Convergence is estimated
using the relative $L^{2}$ error $\epsilon$ on macroscopic variables
with respect to the analytical solution at $t_{max}\approx0.12$.
The reference values $(\Delta t_{ref},\epsilon_{ref})$ of the logarithmic
scale are $\Delta t_{ref}=0.024,\epsilon_{ref}=0.0143$.}
\end{figure}

We compared several $1^{st}$ and $2^{nd}$ order splitting schemes
built up from either fully implicit ($\theta=1$) first order or Crank-Nicolson
($\theta=\frac{1}{2}$) steps:
\begin{itemize}
\item Lie first order splitting scheme $M_{1}^{s}=T_{1}(\Delta t)$$S_{1}(\Delta t)$$C_{1}(\Delta t)$
with first order building blocks.
\item Lie first order splitting scheme $M_{1,2}^{s}=T_{2}(\Delta t)$$S_{2}(\Delta t)$$C_{2}(\Delta t)$
with second order building blocks, for which the order loss comes
from the splitting error.
\item a palindromic second order Strang scheme $M_{2}^{s}=T_{2}(\frac{\Delta t}{2})$$S_{2}(\frac{\Delta t}{2})$$C_{2}(\Delta t)$$S_{2}(\frac{\Delta t}{2})$$T_{2}(\frac{\Delta t}{2})$
. 
\item a collapsed version of the second order Strang scheme $M_{2}^{s}$for
which the last transport substep of each global step and the first
transport substep of the next one are fused in a single $T_{2}(\Delta t)$
substep except obviously for the first and last time steps of the
simulation. 
\end{itemize}
We performed time order convergence tests on a $2D$ square mesh partitioned
into $4\times4$ macrocells with $1024$ points per macrocell. As
shown by Fig. \ref{fig:euler_gravity_order}, we obtain the expected
convergence orders for each of the splitting schemes used.  

\subsection{2D Flow around a cylinder using a penalization method}

In this test case we considered the flow of a fluid in a rectangular
duct with a cylindrical solid obstacle. The simulation domain is the
rectangle $[-1,1]\times[-5,5]$ 

\begin{figure}
\begin{centering}
\includegraphics[width=0.9\textwidth]{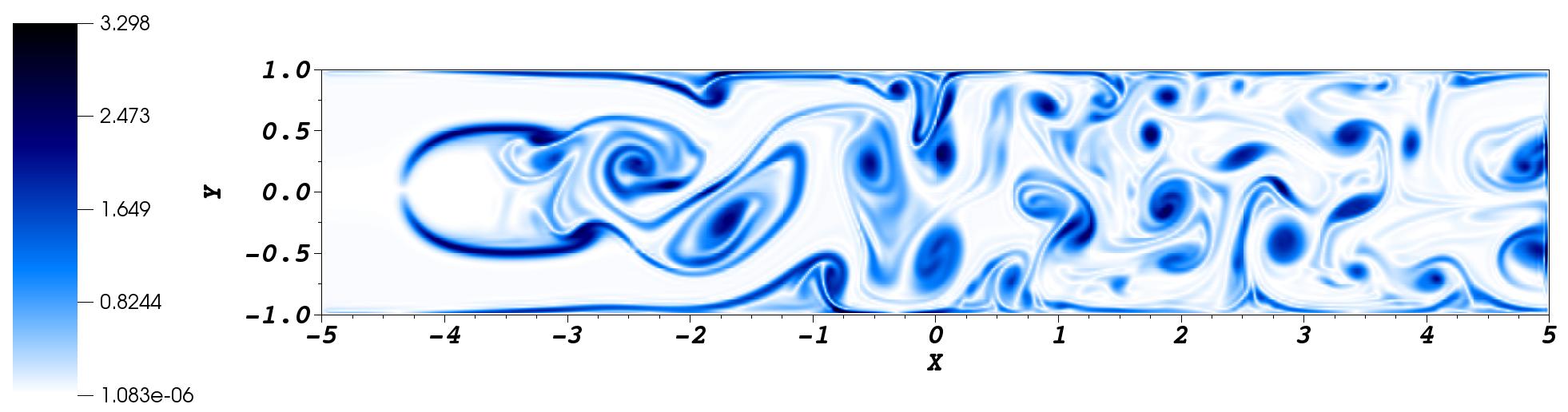}
\par\end{centering}
\caption{\label{fig:euler_cylinder_flow}. Flow around cylindrical obstacle.
Vorticity norm $\vert\nabla\times u\vert$ at $t=83$, showing the
turbulent field behind the obstacle. }
\end{figure}

The effect of the obstacle on the flow is modeled using a volumic
source term of the form 

\begin{equation}
\vs=K(\mathbf{x})(\vw-\vw_{s}),
\end{equation}

with $\vw_{s}=[1.0,0,0]^{t}$ the target fluid state in the ``solid''
part of the domain and the relaxation frequency $K(\mathbf{x})$ is
given by

\begin{equation}
K(\vx)=K_{s}\exp(-\kappa(\vx-\vx_{c})^{2}),
\end{equation}
with $\text{\ensuremath{K_{s}=300}, }\vx_{c}=[-4,0]^{t}$and $\kappa=40$.
The net effect is a very stiff relaxation towards a flow with zero
velocity and the reference density near the center $\vx_{c}$ of the
frequency mask. The effective diameter of the cylinder for this simulation
is about $0.5$. The initial state, which is also applied at the duct
boundaries for the whole simulation is $\rho=1,u_{x}=0.03,u_{y}=0$.
Accounting for the fact that for this model the sound speed is $1/\sqrt{3}$,
the Mach number of the unperturbed flow is approximately $0,017$.
The simulation was performed on a macromesh with $16\times16$ macrocells
stretched with a $1:5$ aspect ratio to match the domain dimension;
each macrocell contains $12\times60$ integration points. On figure
\ref{fig:euler_cylinder_flow} we show the vorticity norm at $t=83$,
when turbulence is well developed in the wake of the obstacle. 

\section{Conclusion}

In this paper, we have presented an optimized implementation of the
Palindromic Discontinuous Galerkin Method for solving kinetic equations
with stiff relaxation. The method presents the following interesting
features:
\begin{itemize}
\item It can be used for solving any hyperbolic system of conservation laws.
\item It is asymptotic-preserving with respect to the stiff relaxation.
\item It is implicit and thus is not limited by CFL conditions.
\item Despite being formally implicit, it requires only explicit computations. 
\item It is easy to increase the time order with a composition method.
\item It presents many opportunities for parallelization and optimization:
in this paper we have presented the parallelization of the method
with the aid of the MPI version of the StarPU runtime system. In this
way we address both shared memory and distributed memory MIMD parallelism.
\end{itemize}
Our perspectives are now to apply the method for computing MHD instabilities
in tokamaks. We will also try to extend the method to more general
boundary conditions.

\bibliographystyle{plain}
\bibliography{schnaps_cemracs}

\end{document}